\documentclass{amsart}

\makeatletter
\def\subsection{\@startsection{subsection}{2}%
  \z@{.5\linespacing\@plus.7\linespacing}{.5\linespacing}%
  {\normalfont\bfseries}}
\makeatother

\makeatletter
\def\@defaultbiblabelstyle#1{[#1]}
\makeatother

\makeatletter
\def\@setauthors{%
  \begingroup
  \def\thanks{\protect\thanks@warning}%
  \trivlist
  \centering\footnotesize \@topsep30\p@\relax
  \advance\@topsep by -\baselineskip
  \item\relax
  \author@andify\authors
  \def\\{\protect\linebreak}%
  \authors%
  \ifx\@empty\contribs
  \else
    ,\penalty-3 \space \@setcontribs
    \@closetoccontribs
  \fi
  \endtrivlist
  \endgroup
}
\def\@settitle{\begin{center}%
  \baselineskip14\p@\relax
    \bfseries
  \@title
  \end{center}%
}
\makeatother

\usepackage{amssymb}
\usepackage[colorlinks=true, citecolor=blue]{hyperref}
\usepackage{upref}
\usepackage{dsfont}
\usepackage{upgreek}
\usepackage{extarrows}
\usepackage[all]{xy}
\allowdisplaybreaks[4]
\usepackage{caption}
\usepackage{tikz}
\usetikzlibrary{calc, graphs, positioning, quotes}

\newtheorem{theorem}{Theorem}[section]
\newtheorem{lemma}[theorem]{Lemma}

\theoremstyle{definition}

\theoremstyle{remark}
\newtheorem{remark}[theorem]{Remark}

\numberwithin{equation}{section}
\allowdisplaybreaks[4]

\begin{document}

\title[Some infinite-dimensional representations]{Some infinite-dimensional representations of certain Coxeter groups}

\author[Hongsheng Hu]{Hongsheng Hu}
\address{Beijing International Center for Mathematical Research, Peking University, No.\ 5 Yiheyuan Road, Haidian District, Beijing 100871, China}
\email{huhongsheng(at)amss(dot)ac(dot)cn}
\urladdr{\href{https://huhongsheng.github.io/}{https://huhongsheng.github.io}}




\date{March 23, 2024}


\begin{abstract}
A Coxeter group admits infinite-dimensional irreducible complex representations if and only if it is not finite or affine.
In this paper, we provide a construction of some of those representations for certain Coxeter groups using some topological information of the corresponding Coxeter graphs.
\end{abstract}

\maketitle

\setcounter{tocdepth}{2}
\tableofcontents

\section{Introduction}
  Let $(W,S)$ be an irreducible Coxeter group of finite rank, i.e., its Coxeter graph is connected and $\lvert S \rvert < \infty$.
  If $W$ is finite, then irreducible representations (over $\mathbb{C}$) of $W$ are certainly finite dimensional.
  If $W$ is an affine Weyl group, then it is also well known that its irreducible representations are of finite dimension (one may refer to \cite{Kato83}, \cite[proof of Prop.\ 5.13]{KL87}, \cite[Prop.\ 1.2]{Xi07} for more details).
  In general, we have the following fact.

  \begin{theorem}\label{thm1.1}
     All irreducible complex representations of $W$ are of finite dimension if and only if $W$ is a finite group or an affine Weyl group.
  \end{theorem}

  The author owes a proof of this theorem to an anonymous referee of a previous version of this paper (see the appendix).
  Nevertheless, the proof only tells us the existence of infinite-dimensional irreducible representations of infinite non-affine Coxeter groups, but it fails to construct such representations.

  The main aim of this paper is to construct some irreducible representations of infinite dimension of a Coxeter group $(W,S)$ satisfying either of the following:
  \begin{enumerate}
      \item \label{case1} there are at least two circuits in the Coxeter graph;
      \item \label{case2} there is at least one circuit in the Coxeter graph, and $m_{st} \ge 4$ for some $s,t \in S$ (for  $s,t \in S$, we denote by $m_{st}$ the order of $st$).
  \end{enumerate}
  The main idea is to glue together many copies of representations of different dihedral subgroups of $W$, so that they form a ``big'' representation of $W$.
  The way of gluing is encoded in some topological information of the Coxeter graph.
  This method is inspired by the author's previous work \cite{Hu21}.

  The paper is organized as follows.
  Section \ref{sec-prel} records some basic facts about representations of dihedral groups, as well as coverings and fundamental groups of graphs.
  Section \ref{sec-3} deals with case \ref{case1}, in which the fundamental group of the Coxeter graph is a non-abelian free group.
  We utilize an infinite-dimensional irreducible representation of this free group to do the ``gluing''.
  In Section \ref{sec-4}, we use the universal covering of the Coxeter graph to achieve our goal for case \ref{case2}.
  In Section \ref{sec-5}, we give another example of infinite-dimensional irreducible representation of a specific Coxeter group whose Coxeter graph has no circuits.
  Finally, in the appendix we present the proof of Theorem \ref{thm1.1} which is given by an anonymous referee.

\section{Preliminaries}\label{sec-prel}
  In this section we recollect some notations and terminology used in this paper.
  We use $e$ to denote the identity in a group.
  Coxeter groups considered throughout this paper are all irreducible and of finite rank.

\subsection{Representations of dihedral groups}\label{subsec2.1}
  For a finite dihedral group $D_m := \langle r,t \mid r^2 = t^2 = (rt)^m = e \rangle$, we denote by $\mathds{1}$ and $\varepsilon$, respectively, the trivial and the sign representation, i.e., $\mathds{1} : r,t \mapsto 1$, $\varepsilon : r,t \mapsto -1$.
  If $m$ is even, there are two more representations of dimension 1, i.e.,
  \begin{equation*}
    \varepsilon_r: r \mapsto -1, t \mapsto 1; \quad \quad \varepsilon_t: r \mapsto 1, t \mapsto -1.
  \end{equation*}

  Let $\mathbb{C} \beta_r \oplus \mathbb{C} \beta_t$ be a vector space with formal basis $\{\beta_r, \beta_t\}$.
  For any integer $k$ satisfying $1 \le k < m/2$, let $\rho_k$ denote the irreducible representation of $D_m$ on $\mathbb{C} \beta_r \oplus \mathbb{C} \beta_t$ defined by
  \begin{align*}
    r \cdot \beta_r & = - \beta_r, & r \cdot \beta_t & = \beta_t + 2\cos\frac{k\uppi}{m} \beta_r,\\
    t\cdot \beta_t & = - \beta_t, & t \cdot \beta_r &= \beta_r + 2\cos\frac{k\uppi}{m} \beta_t.
  \end{align*}
  Intuitively, $r$ and $t$ act on the (real) plane by two reflections with respect to two lines with an angle of $\frac{k \uppi}{m}$; see Figure \ref{rho-k}.

\begin{figure}[ht]
    \centering
    \begin{tikzpicture}
      \draw[dashed] (-1.3,0)--(2,0);
      \draw[dashed] (240:1.3)--(60:2);
      \draw[->] (0,0)--(0,1.5);
      \draw[->] (0,0)--(330:1.5);
      \draw (0.3,0) arc (0:60:0.3);
      \draw[<->] ($(1.9,0) + (300:0.2)$) arc (300:420:0.2);
      \draw[<->] ($(60:1.9) + (0.2,0)$) arc (0:120:0.2);
      \node[left] (ar) at (0,1.5) {$\beta_r$};
      \node[below] (at) at (330:1.5) {$\beta_t$};
      \node[right] (r) at (2.1,0) {$r$};
      \node[right] (t) at (63:2.3) {$t$};
      \node[right] (angle) at (0.2,0.3) {$\frac{k\uppi}{m}$};
    \end{tikzpicture}
    \caption{The representation $\rho_k: D_m \to \operatorname{GL}(\mathbb{C} \beta_r \oplus \mathbb{C} \beta_t)$}\label{rho-k}
\end{figure}

  \begin{remark}\label{rmk2.2}
    If $m$ is even and $k = m/2$, we may define $\rho_{m/2}$ as well by the same formulas, but then $\rho_{m/2} \simeq \varepsilon_r \oplus \varepsilon_t$ is reducible.
  \end{remark}

  \begin{remark}
    We have  described the full set of irreducible representations of $D_m$, namely,
    \begin{gather*}
      \{\mathds{1}, \varepsilon\} \cup \{\rho_1, \dots, \rho_{\frac{m-1}{2}}\}, \quad \text{if $m$ is odd}; \\
      \{\mathds{1}, \varepsilon, \varepsilon_r, \varepsilon_t\} \cup \{\rho_1, \dots, \rho_{\frac{m}{2}-1}\}, \quad \text{if $m$ is even}.
    \end{gather*}
  \end{remark}

  The following lemma will be frequently used in Sections \ref{sec-3} and \ref{sec-4} to prove our ``gluing'' is feasible.

  \begin{lemma}\label{2.2}
    The $+1$-eigenspaces of $r$ and $t$ in $\rho_k$ are both one dimensional.
    However, there is no nonzero vector that can be fixed by $r$ and $t$ simultaneously.
  \end{lemma}

\subsection{Graphs and the universal covering}
  The way we glue these representations will be encoded in some topological information of the Coxeter graph.

  By definition, an \emph{(undirected) graph} $G = (S,E)$ consists of a set $S$ of vertices and a set $E$ of edges, and elements in $E$ are of the form $\{s,t\} \subseteq S$ (unordered).
  For our purpose, we only consider graphs without loops and multiple edges, i.e., there is no edge of the form $\{s,s\}$, and each edge $\{s,t\}$ occurs at most once in $E$.
  In a Coxeter graph, $m_{st}$ is regarded as a label on the edge rather than a multiplicity.
  We say $G$ is a \emph{finite graph} if $S$ is a finite set.

  A sequence $(s_1, s_2, \dots, s_n)$ of vertices is called a \emph{path} in $G$ if $\{s_i, s_{i+1}\} \in E, \forall i$.
  If $s_1 = s_n$, then we say such a path is a \emph{closed path}.
  If further $s_1, \dots, s_{n-1}$ are distinct in this closed path, then the path is called a \emph{circuit}.

  If every two vertices can be connected by a path, then we say $G$ is a \emph{connected graph}.
  A connected graph without circuits is called a \emph{tree}.
  For a connected graph $G = (S,E)$, if $T = (S, E_0)$ is a tree with the same vertices set $S$ and $E_0 \subseteq E$, then $T$ is called a \emph{spanning tree} of $G$.
  This condition is equivalent to say $\lvert E_0 \rvert = \lvert S \rvert - 1$ when $G$ is a connected finite graph.
  Any connected graph has a spanning tree, but not unique in general.

  Let $G = (S, E)$ and $G^\prime = (S^\prime, E^\prime)$ be two graphs.
  If $p : S \to S'$ is a map of sets such that for any edge $\{s,t\}$ in $G$ we have $p(s) \ne p(t)$ and $\{p(s), p(t)\}$ is also an edge in $G'$, then we say that $p$ is a \emph{morphism} of graphs.
  We simply denote a morphism by $p : G \to G'$.

  For any $s \in S$, we denote
  \begin{equation*}
    E_s := \{t \in S \mid \{s,t\} \in E\}.
  \end{equation*}
  Suppose $G'$ is connected and finite, and suppose $p : G \to G'$ is a morphism.
  If $p(S) = S'$, and if for any $s \in S$ the restriction of $p$ to $E_s$ gives rise to a bijection $E_s \xrightarrow{\sim} E'_{p(s)}$, then $p$ is called a \emph{covering}.
  It is natural to regard $G^\prime$ and $G$ as locally finite simplicial complexes.
  Then  $p$ is also a covering of topological spaces.

  Conversely, we view $G$ as a topological space and suppose $p: X \to G$ is a covering of the topological space $G$.
  Then $X$ has a graph structure such that $p$ is a morphism of graphs (see \cite[Theorem 83.4]{Munkres00}).
  In particular, if $p$ is the universal covering, then $X$ is a tree.
  Thus, for any connected finite graph, we can talk about its universal covering graph.

\subsection{The fundamental group of a graph} \label{2.3}
  Let $G = (S, E)$ be a connected graph, and $T = (S, E_0)$ be a spanning tree.
  For any edge $\mathfrak{e} \in E \setminus E_0$, if we choose a vertex $s_\mathfrak{e}$ of $\mathfrak{e}$ to be its head and the other $t_\mathfrak{e}$ to be its tail, then $\mathfrak{e} = \{s_\mathfrak{e}, t_\mathfrak{e}\}$ and there is a unique circuit $c^\prime_\mathfrak{e}$ in $(S, E_0 \cup \{\mathfrak{e}\})$ of the form $c^\prime_\mathfrak{e} = (s_\mathfrak{e}, t_\mathfrak{e}, \dots, s_\mathfrak{e})$.


  Fix a vertex $s_0 \in S$.
  For any $c^\prime_\mathfrak{e}$, there is a unique path in $T$ without repetitive vertices from $s_0$ to $s_\mathfrak{e}$.
  We denote the path by $p_\mathfrak{e}$.
  Define $c_\mathfrak{e}$ to be the concatenation of $p_\mathfrak{e}, c^\prime_\mathfrak{e}, p_\mathfrak{e}^{-1}$, where $p_\mathfrak{e}^{-1}$ is the inverse path in the obvious sense.
  Then each $c_\mathfrak{e}$ is a closed path from $s_0$ to itself.
  If we view $G$ as a topological space, then we have the following result on its fundamental group $\pi_1(G)$.
  \begin{lemma} [{\cite[Theorem 84.7]{Munkres00}}] \label{lem-2.5} 
    $\pi_1(G)$ is a free group with a set of free generators $\{c_\mathfrak{e} \mid \mathfrak{e} \in E \setminus E_0\}$.
    In particular, if there is more than one circuit in $G$, then $\pi_1(G)$ is non-abelian.
  \end{lemma}
  \begin{remark}\label{rmk2.6}
    Note that $p_\mathfrak{e}$ is a path in $T$.
    Thus, all edges in $c_\mathfrak{e}$ except $\mathfrak{e}$ lie in $E_0$, and $\mathfrak{e}$ appears in $c_\mathfrak{e}$ only once.
    If $\mathfrak{e}^\prime \in E\setminus E_0$ is another edge, then $\mathfrak{e}^\prime$ does not appear in $c_\mathfrak{e}$.
  \end{remark}

\section{Representations via fundamental groups of Coxeter graphs} \label{sec-3}
  From now on, suppose $(W,S)$ is an irreducible Coxeter group of finite rank with Coxeter graph $G = (S,E)$.
  Then, $G$ is a connected finite graph.

  In this section we assume that
  \begin{equation*}
    \text{there are at least two circuits in $G$.}
  \end{equation*}
  Thus, by Lemma \ref{lem-2.5}, $\pi_1(G)$ is a non-abelian free group.
  In this section we use an infinite-dimensional irreducible representation of $\pi_1(G)$ to construct another such representation for $W$.

  For convenience, we may further assume
  \begin{equation*}
    m_{st} < \infty, \quad \forall s,t \in S.
  \end{equation*}
  This is not essential.
  If some $m_{st}$'s are infinity, then we replace them by any integer larger than 2 (e.g., 3), so that we obtain another Coxeter group $(W_1,S)$ and a surjective homomorphism $W \twoheadrightarrow W_1$, $s \mapsto s$.
  Ignoring labels on edges, the two Coxeter groups have the same Coxeter graph in a topological sense.
  An irreducible representation of $W_1$ becomes an irreducible representation of $W$ via pulling back by the homomorphism.

\subsection{The construction} \label{3.1}
  We fix $s_0 \in S$.
  Let $T = (S, E_0)$ be a spanning tree of $G = (S,E)$.
  Then we obtain a set of free generators of $\pi_1(G)$ by the method in Subsection \ref{2.3}, say $c_1, \dots, c_l$ ($l < \infty$ since $G$ is a finite graph), with $s_0$ on each of them.

  We have $\lvert E \setminus E_0 \rvert \ge 2$ since we assume that there is more than one circuit in $G$.
  We then fix two distinct edges:
  \begin{equation*}
    \{s_1, t_1\}, \{s_2, t_2\} \in E \setminus E_0.
  \end{equation*}
  We may assume that $\{s_1, t_1\}, \{s_2, t_2\}$ lie in $c_1, c_2$ respectively, and $c_1$ goes through $t_1$ first and then $s_1$, and $c_2$ goes through $t_2$ first and then $s_2$.
  In our choice of $c_i$, the edge $\{s_1, t_1\}$ appears in $c_1$ only once, while it does not appear in other $c_i$'s (see Remark \ref{rmk2.6}).
  Similar for $\{s_2, t_2\}$.
  The two edges might share a common vertex, like $s_1 = s_2$, but it does not matter.

  We define a vector space $V$ with formal basis $\{\alpha_{s,n} \mid n \in \mathbb{Z}, s \in S \}$,
  \begin{equation*}
    V := \bigoplus_{n \in \mathbb{Z}, s \in S} \mathbb{C} \alpha_{s,n}.
  \end{equation*}
  For any $s,t \in S$ and $n \in \mathbb{Z}$, we define $s \cdot \alpha_{t,n}$ as follows:
  \begin{enumerate}
    \item if $s = t$, then $s \cdot \alpha_{s,n} : = - \alpha_{s,n}$;
    \item $s_1 \cdot \alpha_{t_1, n} := \alpha_{t_1, n} + 2 \cos \frac{\uppi}{m_{s_1t_1}} \alpha_{s_1, n+1}$,
        \newline
        $t_1 \cdot \alpha_{s_1, n+1} := \alpha_{s_1, n+1} + 2 \cos \frac{\uppi}{m_{s_1t_1}} \alpha_{t_1,n}$;
        \newline
        $s_2 \cdot \alpha_{t_2, n} := \alpha_{t_2, n} + 2^{n+1} \cos \frac{\uppi}{m_{s_2t_2}} \alpha_{s_2, n+1}$,
        \newline
        $t_2 \cdot \alpha_{s_2, n+1} := \alpha_{s_2, n+1} + 2^{-n+1} \cos \frac{\uppi}{m_{s_2t_2}} \alpha_{t_2,n}$;
    \item if the unordered pair $\{s,t\} \ne \{s_1, t_1\}$ or $\{s_2, t_2\}$, and if $s \ne t$, then \newline $s \cdot \alpha_{t,n} := \alpha_{t,n} + 2 \cos \frac{\uppi}{m_{st}} \alpha_{s,n}$.
  \end{enumerate}

  Looking at case (2), one can see that the vectors $\alpha_{t_2,n}$ and $2^n \alpha_{s_2, n+1}$ span an irreducible representation isomorphic to $\rho_1$ (see Subsection \ref{subsec2.1}) of the dihedral subgroup $\langle s_2, t_2 \rangle$.
  The vectors $\alpha_{t_2,n}$ and $2^n \alpha_{s_2, n+1}$ play the same roles as the $\beta$'s in Subsection \ref{subsec2.1}.
  Similarly, $\{\alpha_{t_1,n}, \alpha_{s_1,n+1}\}$ span an irreducible representation isomorphic to $\rho_1$ of $\langle s_1, t_1 \rangle$.

  In case (3), if $m_{st} = 2$, then $s \cdot \alpha_{t,n} = \alpha_{t,n}$; if $m_{st} \ge 3$, then $\{\alpha_{t,n}, \alpha_{s,n}\}$ span an irreducible representation isomorphic to $\rho_1$ of $\langle s,t \rangle$.

  \begin{lemma} \label{3.2}
    $V$ is a representation of $W$ with the action defined above.
  \end{lemma}

  \begin{proof}
    Obviously, $s^2$ acts as the identity for any $s \in S$.

    For $s,t \in S$ and $m_{st} = 2$, we need to show $st \cdot \alpha_{r,n} = ts \cdot \alpha_{r,n}$, $\forall r \in S$.
    Note that we have
    \[s \cdot \alpha_{t, n} = \alpha_{t,n} \text{ and } t \cdot \alpha_{s,n} = \alpha_{s,n} \text{ for any } n \in \mathbb{Z},\]
    since $m_{st} = 2$.
    If $r = s$ or $r = t$, then clearly it holds $st \cdot \alpha_{r,n} = ts \cdot \alpha_{r,n}$.
    If $r \ne s$ and $r \ne t$, then
    \[s \cdot \alpha_{r,n} = \alpha_{r,n} + c_1 \alpha_{s, n_1}
    \text{ and }
    t \cdot \alpha_{r,n} = \alpha_{r,n} + c_2 \alpha_{t, n_2}\]
    for some $c_1, c_2 \in \mathbb{C}$ and $n_1, n_2 \in \{n, n \pm 1\}$, and then
    \begin{align*}
      st \cdot \alpha_{r,n} & = s \cdot (\alpha_{r,n} + c_2 \alpha_{t ,n_2}) = \alpha_{r,n} + c_1 \alpha_{s, n_1} + c_2 \alpha_{t ,n_2}, \\
      ts \cdot \alpha_{r,n} & = t \cdot (\alpha_{r,n} + c_1 \alpha_{s, n_1}) = \alpha_{r,n} + c_2 \alpha_{t, n_2} + c_1 \alpha_{s, n_1}.
    \end{align*}
    Therefore, we have  $st \cdot \alpha_{r,n} = ts \cdot \alpha_{r,n}$ as desired.

    Now assume $m_{st} \ge 3$.
    We need to verify that $(st)^{m_{st}} \cdot \alpha_{r,n} = \alpha_{r,n}$.
    This is also obvious if $r = s$ or $r = t$, since we are in the dihedral world.
    If $s, t, r$ are distinct, then we have the following cases.

    (1) If any two of $s,t,r$ are not $\{s_1, t_1\}$ or $\{s_2, t_2\}$, then the three-dimensional subspace spanned by $\alpha_{r,n}$, $\alpha_{s,n}$, $\alpha_{t,n}$, denoted by $U$, stays invariant under the actions of $s$ and $t$.
    We write $U_s := \{v \in U \mid s \cdot v = v\}$, $U_t := \{v \in U \mid t \cdot v = v\}$.
    Then $\dim U_s = \dim U_t = 2$, and thus there exists $0 \ne v_0 \in U$ such that $s \cdot v_0 = t \cdot v_0 = v_0$.
    Note that $3 \le m_{st} < \infty$ and $\mathbb{C} \alpha_{s,n} \oplus \mathbb{C} \alpha_{t,n}$ forms a representation isomorphic to $\rho_1$ of $\langle s,t \rangle$.
    Hence, $v_0 \notin \mathbb{C} \alpha_{s,n} \oplus \mathbb{C} \alpha_{t,n}$ by Lemma \ref{2.2}, and then $\{ v_0, \alpha_{s,n}, \alpha_{t,n}\}$ is a basis of $U$.
    Now we can see that $(st)^{m_{st}} \cdot \alpha_{r,n} = \alpha_{r,n}$.

    (2) If $m_{rt} = 2$, then there exists $k \in \{n-1, n, n+1\}$ and $q \in \{k-1, k, k+1\}$ such that $\alpha_{r,n}, \alpha_{s,k}, \alpha_{t,q}$ span an $s,t$-invariant subspace.
    By the same arguments in case (1), we have $(st)^{m_{st}} \cdot \alpha_{r,n} = \alpha_{r,n}$.
    The case $m_{rs} = 2$ is similar.

    In the following cases, we assume $s,t,r$ do not commute with each other.

    (3) If $s = s_1$, $t = t_1$, while $s_2, t_2$ do not occur simultaneously in $s,r,t$, then $\alpha_{s_1, n+1}$, $\alpha_{t_1,n}$, $\alpha_{r,n}$, $\alpha_{s_1, n}$, $\alpha_{t_1, n-1}$ span a five-dimensional $s,t$-invariant subspace $U$.
    Define $U_s$, $U_t$ as in case (1), then $\dim U_s = \dim U_t = 3$.
    The same argument shows that $(s_1t_1)^{m_{s_1t_1}} \cdot \alpha_{r,n} = \alpha_{r,n}$.

    (4) If $s = s_1$, $r = t_1$, while $s_2, t_2$ do not occur simultaneously in $s,r,t$, then $\alpha_{s_1,n}$, $\alpha_{t,n}$, $\alpha_{r,n}$, $\alpha_{s_1, n+1}$, $\alpha_{t, n+1}$ span an $s,t$-invariant subspace.
    The same argument works.

    (5) If $s = s_1$, $t = t_1 = t_2$, $r = s_2$, then $\alpha_{r,n}, \alpha_{s,n}, \alpha_{t, n-1}$ span an $s,t$-invariant subspace.
    The same argument works.

    (6) If $s = s_1$, $t = t_1 = s_2$, $r = t_2$, then $\alpha_{t_1, n-1}$, $\alpha_{s_1,n}$, $\alpha_{r,n}$, $\alpha_{s_2, n+1}$, $\alpha_{s_1, n+2}$ span an $s,t$-invariant subspace.
    The same argument works.

    (7) If $s = s_1$, $t = s_2$, $r = t_1 = t_2$, then $\alpha_{s_1, n+1}$, $\alpha_{r,n}$, $\alpha_{s_2, n+1}$ span an $s,t$-invariant subspace.
    The same argument works.

    (8) If $s = s_1$, $t = t_2$, $r = t_1 = s_2$, then $\alpha_{t_2,n+1}$, $\alpha_{s_1, n+1}$, $\alpha_{r,n}$, $\alpha_{t_2, n-1}$, $\alpha_{s_1, n-1}$ span an $s,t$-invariant subspace.
    The same argument works.

    In the above cases, if we exchange the letters $s,t$ or the indices $1,2$, then the arguments are totally the same.
    Thus, we always have $(st)^{m_{st}} \cdot \alpha_{r,n} = \alpha_{r,n}$.
  \end{proof}

\subsection{The infinite-dimensional irreducible quotient}

  \begin{theorem}\label{thm3.2}
    Recall that $(W,S)$ is irreducible, and that there are at least two circuits in its Coxeter graph $G$.
    Let $V$ be defined as in Subsection \ref{3.1}, and
    \begin{equation*}
      V_0 := \{v \in V \mid s \cdot v = v, \forall s \in S\}.
    \end{equation*}
    Then the representation $V / V_0$ of $W$ is irreducible of infinite dimension.
  \end{theorem}

  \begin{proof}
    We denote $V_-^s := \{v \in V \mid s \cdot v = -v\}$.
    Then $V_-^s = \bigoplus_n \mathbb{C} \alpha_{s,n}$.
    For any edge $\{s,t\}$ in $G$, $V_-^s \oplus V_-^t$ is a subrepresentation of $\langle s,t \rangle$ in $V$, isomorphic to an infinite direct sum of $\rho_1$.
    For any $v \in V_-^s$, let $f_{st}(v) := (t \cdot v - v) / 2 \cos \frac{\uppi}{m_{st}}$.
    Then $f_{st}(v) \in V_-^t$, and the linear map $f_{st}: V_-^s \to V_-^t$ is a linear isomorphism of vector spaces.
    For example, when  $\{s,t\} \ne \{s_1,t_1\}$ or $\{s_2,t_2\}$, we have $f_{st}(\alpha_{s,n}) = \alpha_{t,n}$.
    Moreover, $f_{st}(v)$ lies in the subrepresentation generated by $v$.

    Let $0 \ne v \in V$, and let $U$ be the subrepresentation generated by $v$.
    If $v \notin V_0$, say, $t \cdot v \ne v$, then $t \cdot v - v \in V_-^t \cap U$.
    Suppose $(r_0 = t, r_1, \dots, r_k = s_0)$ is a path connecting $t$ and $s_0$.
    Here $s_0$ is the vertex fixed in Subsection \ref{3.1}.
    Then,
    \begin{equation*}
      v_0 : = f_{r_{k-1}r_k} \cdots f_{r_1r_2} f_{r_0 r_1} (t \cdot v - v) \in V_-^{s_0} \cap U \text{ and } v_0 \ne 0.
    \end{equation*}

    Apply the maps $f_{**}$ along the closed paths $c_1, \dots, c_l$ chosen in Subsection \ref{3.1}.
    Then we obtain $l$ linear isomorphisms of $V_-^{s_0}$, denoted by $X_1, \dots, X_l$, respectively.
    This makes $V_-^{s_0}$ form a representation of the free group $\pi_1(G)$.
    Except $X_1$ and $X_2$, other $X_i$'s are identity maps of $V_-^{s_0}$, and we have
    \begin{equation*}
      X_1 (\alpha_{s_0, n}) = \alpha_{s_0, n+1}, \quad X_2 (\alpha_{s_0, n}) = 2^n \alpha_{s_0, n+1}.
    \end{equation*}
    It is easy to verify that $V_-^{s_0}$ is an irreducible representation of $\pi_1(G)$.
    Notice that $0 \ne v_0 \in V_-^{s_0} \cap U$, $X_i^{\pm 1} (v_0) \in U$.
    Thus, $V_-^{s_0} \subseteq U$.
    Since $G$ is connected, $V_-^{s_0}$ generates the whole representation $V$.
    Hence, $V = U$.
    We have proved that $V / V_0$ is an irreducible representation of $W$.

    Note that $s_0$ acts on $V_-^{s_0}$ by $-1$.
    So $V_-^{s_0} \cap V_0 = 0$, and thus $\dim V/V_0 = \infty$.
  \end{proof}

\section{Representations via universal coverings of Coxeter graphs} \label{sec-4}
  In this section, we assume that
  \begin{gather*}
    \text{the Coxeter graph $G = (S,E)$ is not a tree,} \\
    \text{and there exist $s_1, s_2 \in S$ such that $m_{s_1s_2} \ge 4$.}
  \end{gather*}
  In this section we use the universal covering of $G$ to construct an infinite-dimensional representation of $W$, then find an irreducible (sub)quotient in it.

  For the same reason stated before Subsection \ref{3.1}, we may further assume
   \begin{equation*}
     m_{st} < \infty, \quad \forall s,t \in S.
   \end{equation*}

\subsection{The construction} \label{4.1}
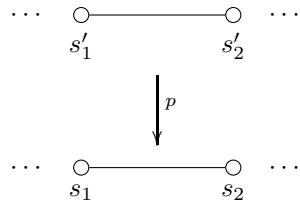
  \begin{figure} [ht]
    \centering
    \begin{tikzpicture}
      \node [circle, draw, inner sep=2pt, label=below:$s_2^\prime$] (s1p) at (1,2) {};
      \node [circle, draw, inner sep=2pt, label=below:$s_1^\prime$] (s0p) at (-1,2) {};
      \node (lp) at (-1.7,2) {$\cdots$};
      \node (rp) at (1.7,2) {$\cdots$};
      \graph{(s1p) -- (s0p);};

      \draw [->] (0,1.3) -- (0,0.5);
      \node [label=right:$p$] (p) at (-0.1,0.9) {};

      \node [circle, draw, inner sep=2pt, label=below:$s_2$] (s1) at (1,0) {};
      \node [circle, draw, inner sep=2pt, label=below:$s_1$] (s0) at (-1,0) {};
      \node (l) at (-1.7,0) {$\cdots$};
      \node (r) at (1.7,0) {$\cdots$};
      \graph{(s1) -- (s0);};
    \end{tikzpicture}
    \caption{The edge $\{s_1^\prime, s_2^\prime\}$}\label{fig3}
  \end{figure}

  We fix $s_1, s_2 \in S$ such that $m_{s_1s_2} \ge 4$.
  Let $p: G^\prime \to G$ be the universal covering of $G$, where $G^\prime = (S^\prime, E^\prime)$.
  Fix an edge $\{s_1^\prime, s_2^\prime\}$ in $G^\prime$ such that $p(s_1^\prime) = s_1$, $p(s_2^\prime) = s_2$, as shown in Figure \ref{fig3}.

  We define a vector space $V$ with formal basis $\{\alpha_a \mid a\in S^\prime\}$,
  \begin{equation*}
    V := \bigoplus_{a \in S^\prime} \mathbb{C} \alpha_a.
  \end{equation*}
  For any $s \in S$ and $a \in S'$, we define $s \cdot \alpha_a$ as follows:
  \begin{enumerate}
    \item if $s = p(a)$, then $s \cdot \alpha_a : = - \alpha_a$;
    \item if $s = s_1$, $a = s_2^\prime$, then $s_1 \cdot \alpha_{s_2^\prime} : = \alpha_{s_2^\prime} + 2 \cos \frac{2 \uppi}{m_{s_1s_2}} \alpha_{s_1^\prime}$;
    \item if $s = s_2$, $a = s_1^\prime$, then $s_2 \cdot \alpha_{s_1^\prime} : = \alpha_{s_1^\prime} + 2 \cos \frac{2 \uppi}{m_{s_1s_2}} \alpha_{s_2^\prime}$;
    \item if it is not in the cases above, and if $s$ is not adjacent to  $p(a)$ in $G$, then $s \cdot \alpha_a := \alpha_a$;
    \item if it is not in the cases above, and if $s$ is adjacent to $p(a)$ in $G$, then we denote by $b$ the vertex adjacent to $a$ in $p^{-1} (s)$ (see Figure \ref{fig4}), and $s \cdot \alpha_a := \alpha_a + 2 \cos \frac{\uppi}{m} \alpha_b$, where $m := m_{s,p(a)} \ge 3$.

        
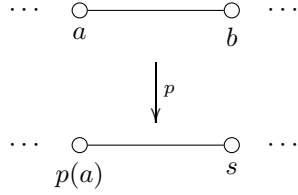
\begin{figure} [ht]
          \centering
          \begin{tikzpicture}
            \node [circle, draw, inner sep=2pt, label=below:$b$] (s1p) at (1,2) {};
            \node [circle, draw, inner sep=2pt, label=below:$a$] (s0p) at (-1,2) {};
            \node (lp) at (-1.7,2) {$\cdots$};
            \node (rp) at (1.7,2) {$\cdots$};
            \graph{(s1p) -- (s0p);};

            \draw [->] (0,1.3) -- (0,0.5);
            \node [label=right:$p$] (p) at (-0.1,0.9) {};

            \node [circle, draw, inner sep=2pt, label=below:$s$] (s1) at (1,0) {};
            \node [circle, draw, inner sep=2pt, label=below:$p(a)$] (s0) at (-1,0) {};
            \node (l) at (-1.7,0) {$\cdots$};
            \node (r) at (1.7,0) {$\cdots$};
            \graph{(s1) -- (s0);};
          \end{tikzpicture}
          \caption{The vertex $b$ in $p^{-1}(s)$ adjacent to $a$}\label{fig4}
        \end{figure}
  \end{enumerate}
  In particular, $\mathbb{C} \alpha_{s_1^\prime} \oplus \mathbb{C} \alpha_{s_2^\prime}$ forms a representation of the dihedral subgroup $\langle s_1,s_2 \rangle$, isomorphic to $\rho_2$ (note that if $m_{s_1s_2} = 4$, then this representation splits, see Remark \ref{rmk2.2}).
  While for other pairs of adjacent vertices $\{a,b\}$ in $G^\prime$, $\mathbb{C} \alpha_a \oplus \mathbb{C} \alpha_b$ forms an irreducible representation of $\langle p(a), p(b) \rangle$ isomorphic to $\rho_1$.

  \begin{lemma} \label{lem-3.3}
    $V$ is a representation of $W$ with the action defined above.
  \end{lemma}

  \begin{proof}
    From the construction, it is clear that $s^2$ acts by identity for any $s \in S$.
    Suppose $s,t \in S$ and $s\ne t$.
    We need to verify that $(st)^{m_{st}} \cdot \alpha_a = \alpha_a$ for any $a \in S^\prime$.
    If $p(a) = s$ or $t$, then $\alpha_a$ lies in a subrepresentation of $\langle s,t \rangle$.
    Thus, we have $(st)^{m_{st}} \cdot \alpha_a = \alpha_a$.
    If $p(a) \neq s$ and $p(a) \ne t$, then the relationship of the three vertices $p(a), s, t$ in $G$ is in one of the following cases (ignoring labels like $m_{st}$ on edges),
    
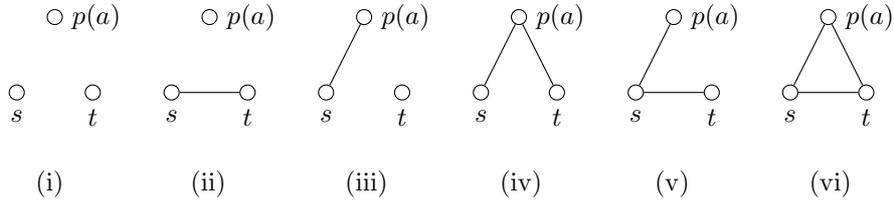
\begin{equation*}
       \begin{tikzpicture}
        \node [circle, draw, inner sep=2pt, label=below:$s$] (s) at (0,0) {};
        \node [circle, draw, inner sep=2pt, label=below:$t$] (t) at (1,0) {};
        \node [circle, draw, inner sep=2pt, label=right:$p(a)$] (pa) at (0.5,1) {};
        \node [label=below:{(i)}] (l1) at (0.5, -0.5) {};
      \end{tikzpicture}  \quad
       \begin{tikzpicture}
        \node [circle, draw, inner sep=2pt, label=below:$s$] (s) at (0,0) {};
        \node [circle, draw, inner sep=2pt, label=below:$t$] (t) at (1,0) {};
        \node [circle, draw, inner sep=2pt, label=right:$p(a)$] (pa) at (0.5,1) {};
        \graph {(s) -- (t);};
        \node [label=below:{(ii)}] (l1) at (0.5, -0.5) {};
      \end{tikzpicture}  \quad
      \begin{tikzpicture}
        \node [circle, draw, inner sep=2pt, label=below:$s$] (s) at (0,0) {};
        \node [circle, draw, inner sep=2pt, label=below:$t$] (t) at (1,0) {};
        \node [circle, draw, inner sep=2pt, label=right:$p(a)$] (pa) at (0.5,1) {};
        \graph {(pa) -- (s);};
        \node [label=below:{(iii)}] (l1) at (0.5, -0.5) {};
      \end{tikzpicture}  \quad
      \begin{tikzpicture}
        \node [circle, draw, inner sep=2pt, label=below:$s$] (s) at (0,0) {};
        \node [circle, draw, inner sep=2pt, label=below:$t$] (t) at (1,0) {};
        \node [circle, draw, inner sep=2pt, label=right:$p(a)$] (pa) at (0.5,1) {};
        \graph {(s) -- (pa) -- (t);};
        \node [label=below:{(iv)}] (l1) at (0.5, -0.5) {};
      \end{tikzpicture}  \quad
      \begin{tikzpicture}
        \node [circle, draw, inner sep=2pt, label=below:$s$] (s) at (0,0) {};
        \node [circle, draw, inner sep=2pt, label=below:$t$] (t) at (1,0) {};
        \node [circle, draw, inner sep=2pt, label=right:$p(a)$] (pa) at (0.5,1) {};
        \graph {(pa) -- (s) -- (t);};
        \node [label=below:{(v)}] (l1) at (0.5, -0.5) {};
      \end{tikzpicture}  \quad
      \begin{tikzpicture}
        \node [circle, draw, inner sep=2pt, label=below:$s$] (s) at (0,0) {};
        \node [circle, draw, inner sep=2pt, label=below:$t$] (t) at (1,0) {};
        \node [circle, draw, inner sep=2pt, label=right:$p(a)$] (pa) at (0.5,1) {};
        \graph {(pa) -- (s) -- (t) -- (pa);};
        \node [label=below:{(vi)}] (l1) at (0.5, -0.5) {};
      \end{tikzpicture}
    \end{equation*}
    (In cases (iii) and (v), exchanging letters $s$ and $t$ does not cause essential differences. So we omit them.)
    In cases (i) (iii) (iv), we may verify directly by definition that $st \cdot \alpha_a = ts \cdot \alpha_a$.
    In case (ii), it is also clear that $(st)^{m_{st}} \cdot \alpha_a = \alpha_a$.

    Suppose we are in case (v).
    We denote by $s^\prime$ the vertex adjacent to $a$ in $p^{-1}(s)$, and by $t^\prime$ the vertex adjacent to $s^\prime$ in $p^{-1}(t)$.
    Then $a$ is not adjacent to $t^\prime$, as shown in Figure \ref{fig5}.
    The three-dimensional subspace spanned by $\alpha_a, \alpha_{s^\prime}, \alpha_{t^\prime}$ stays invariant under the actions of $s$ and $t$.
    By the same arguments as in the proof of Lemma \ref{3.2}, it holds that $(st)^{m_{st}} \cdot \alpha_a = \alpha_a$.

    \begin{figure} [ht]
      \centering
      \begin{tikzpicture}
        \node [circle, draw, inner sep=2pt, label=below:$s^\prime$] (sp) at (0,2) {};
        \node [circle, draw, inner sep=2pt, label=below:$a$] (pap) at (-1.5,2) {};
        \node [circle, draw, inner sep=2pt, label=below:$t^\prime$] (tp) at (1.5,2) {};
        \graph {(pap) -- (sp) -- (tp);};
        \node (lp) at (-2.2,2) {$\cdots$};
        \node (rp) at (2.2,2) {$\cdots$};

        \draw [->] (0,1.3) -- (0,0.5);
        \node [label=right:$p$] (p) at (-0.1,0.9) {};

        \node [circle, draw, inner sep=2pt, label=below:$s$] (s) at (0,0) {};
        \node [circle, draw, inner sep=2pt, label=below:$p(a)$] (pa) at (-1.5,0) {};
        \node [circle, draw, inner sep=2pt, label=below:$t$] (t) at (1.5,0) {};
        \graph {(pa) -- (s) -- (t);};
        \node (l) at (-2.2,0) {$\cdots$};
        \node (r) at (2.2,0) {$\cdots$};
      \end{tikzpicture}
      \caption{The vertices $s^\prime$ and $t^\prime$}\label{fig5}
    \end{figure}

        \begin{figure} [ht]
      \centering
      \begin{tikzpicture}
        \node [circle, draw, inner sep=2pt, label=357:$a$] (a) at (0,0) {};
        \node [circle, draw, inner sep=2pt, label=left:$s^\prime$] (sp) at (-1.5,-0.5) {};
        \node [circle, draw, inner sep=2pt, label=right:$t^\prime$] (tp) at (1.5,-0.5) {};
        \node [circle, draw, inner sep=2pt, label=left:$s^{\prime\prime}$] (spp) at (-1.5,0.5) {};
        \node [circle, draw, inner sep=2pt, label=right:$t^{\prime\prime}$] (tpp) at (1.5,0.5) {};
        \node (d) at (2.5,0) {$\cdots$};
        \node (dd) at (-2.5,0) {$\cdots$};
        \graph {(spp) -- (tpp) -- (a) -- (sp) -- (tp);};

        \draw [->] (0,-1) -- (0, -1.75);
        \node [label=right:$p$] (p) at (-0.1,-1.35) {};

        \node [circle, draw, inner sep=2pt, label=left:$s$] (s1) at (-1.5,-3) {};
        \node [circle, draw, inner sep=2pt, label=right:$t$] (t1) at (1.5,-3) {};
        \node [circle, draw, inner sep=2pt, label=3:$p(a)$] (pa1) at (0,-2.5) {};
        \graph {(pa1) -- (s1) -- (t1) -- (pa1);};
        \node (d1) at (2.5,-2.75) {$\cdots$};
        \node (dd1) at (-2.5,-2.75) {$\cdots$};
      \end{tikzpicture}
      \caption{The vertices $s^\prime, t^\prime, s^{\prime\prime}$ and $t^{\prime\prime}$}\label{fig6}
    \end{figure}

    Suppose we are in case (vi).
    We denote by $s^\prime$ the vertex adjacent to $a$ in $p^{-1}(s)$, by $t^\prime$ the vertex adjacent to $s^\prime$ in $p^{-1}(t)$, by $t^{\prime\prime}$ the vertex adjacent to $a$ in $p^{-1}(t)$, and by $s^{\prime\prime}$ the vertex adjacent to $t^{\prime\prime}$ in $p^{-1}(s)$ (see Figure \ref{fig6}).
    Then $a, t^\prime, s^{\prime\prime}$ are not adjacent to each other.
    Then $\alpha_a, \alpha_{s^\prime}, \alpha_{s^{\prime\prime}}, \alpha_{t^\prime}, \alpha_{t^{\prime\prime}}$ span an $s,t$-invariant subspace of dimension 5.
    The same arguments as in the proof of Lemma \ref{3.2} yield $(st)^{m_{st}} \cdot \alpha_a = \alpha_a$.
  \end{proof}

\subsection{The infinite-dimensional irreducible (sub)quotient} \label{4.3}
  Let $a$ and $b$ be two arbitrary vertices in $G^\prime$.
  Since $G^\prime$ is a tree, there is a unique path $(a = t_0, t_1, \dots, t_n = b)$ connecting $a,b$ such that all $t_i$'s are distinct.
  We define $\operatorname{d}(a,b) := n$ to be the distance between $a$ and $b$.
  We define also
  \begin{equation*}
    S_1^\prime := \{a \in S^\prime \mid \operatorname{d}(a, s_1^\prime) < \operatorname{d}(a, s_2^\prime)\}, \quad S_2^\prime := \{a \in S^\prime \mid \operatorname{d}(a, s_2^\prime) < \operatorname{d}(a, s_1^\prime)\}.
  \end{equation*}
  Then one of them is an infinite set.
  Without loss of generality, we assume $\lvert S_1^\prime \rvert = \infty$.
  Let
  \begin{equation*}
    V_1 := \bigoplus_{a\in S_1^\prime} \mathbb{C} \alpha_a, \quad V_0 := \{v \in V \mid s \cdot v = v, \forall s \in S\}.
  \end{equation*}
  Then $\dim V_1 = \infty$.
  If $m_{s_1s_2} = 4$, then $V_1$ is a subrepresentation of $W$ in $V$.

  \begin{lemma} \label{lem-3.4} \leavevmode
    \begin{enumerate}
      \item If $m_{s_1s_2} > 4$ and $v \in V \setminus V_0$, then $V$ is generated by $v$ as a representation of $W$.
      \item If $m_{s_1s_2} = 4$ and $v_1 \in V_1 \setminus V_0$, then $V_1$ is generated by $v_1$ as a representation of $W$.
    \end{enumerate}
  \end{lemma}

  \begin{proof}
    (1). Suppose $s \cdot v \ne v$, $s \in S$.
    Then by definition we know that $v - s\cdot v$ is a finite sum of the form $\sum_{a \in p^{-1}(s)} x_a \alpha_a$, where each $x_a$ is a complex number.
    Let $u_0 := v - s \cdot v$, and $U$ be the subrepresentation generated by $v$.
    Then $u_0 \in U$.
    We take $a_0 \in p^{-1} (s)$ such that $x_{a_0} \ne 0$.
    Suppose the shortest path in $G^\prime$ connecting $a_0$ and $s_1^\prime$ is \begin{equation*}
      (a_0, a_1, \dots, a_n = s_1^\prime).
    \end{equation*}
    Let $t_i := p(a_i) \in S$, $u_i := u_{i-1} - t_{i} \cdot u_{i-1}$.
    Inductively, we can see that $u_i$ is of the form $\sum_{b \in p^{-1}(t_i)} x_{i,b} \alpha_b$, where $x_{i,b} \in \mathbb{C}$.
    Since $p$ is a covering map, there is only one vertex (namely, $a_{i-1}$) in $p^{-1}(t_{i-1})$ adjacent to $a_i$.
    Thus in the expression of $u_i$, the coefficient $x_{i,a_i}$ of $\alpha_{a_i}$ is nonzero.
    In particular, taking $i = n$, we know that $u_n \in U$ and the coefficient of $\alpha_{s_1^\prime}$ is nonzero.

    We view $V$ as a representation of the finite dihedral group $D := \langle s_1, s_2 \rangle$.
    Since the group algebra $\mathbb{C}[D]$ is semisimple, $V$ decomposes into a direct sum of some copies of irreducible representations of $D$.
    From the construction of $V$, the only  irreducible representations of $D$ which are possible to occur in $V$ are $\mathds{1}, \rho_1, \rho_2$.
    Moreover, $\rho_2$ appears only once, namely, $\mathbb{C} \alpha_{s_1^\prime} \oplus \mathbb{C} \alpha_{s_2^\prime}$.
    Therefore, there is an element $d \in \mathbb{C}[D]$ such that $d \cdot u_n = \alpha_{s_1^\prime}$, and hence $\alpha_{s_1^\prime} \in U$.
    Note that $G^\prime$ is a connected graph.
    From the definition of $V$, we know that $\alpha_{s_1^\prime}$ generates the whole $V$.
    Thus, $U = V$.

    (2). The proof is similar.
    As above, we take a vertex $a_0 \in p^{-1}(s) \cap S_1^\prime$ such that $\alpha_{a_0}$ has nonzero coefficient in the linear expression of $v_1 - s \cdot v_1$ ($\ne 0$), and do the same discussion along the shortest path connecting $a_0$ and $s_1^\prime$ (note that all of the vertices in this path belong to $S_1^\prime$).
    Then we know that in the subrepresentation generated by $v_1$, there is a vector $u_n$ with nonzero coefficient of $\alpha_{s_1^\prime}$.

    Decompose $V_1$ into a direct sum of irreducible representations of $D = \langle s_1, s_2 \rangle$.
    Then only $\mathds{1}, \rho_1, \varepsilon_{s_1}$ occur.
    Moreover, the subrepresentation $\varepsilon_{s_1}$, spanned by $\alpha_{s_1^\prime}$, is of multiplicity one.
    Thus, $\alpha_{s_1^\prime}$ lies in the subrepresentation generated by $v_1$, while $\alpha_{s_1^\prime}$ generates the whole representation $V_1$.
  \end{proof}

  Let
  \begin{gather*}
    \overline{V} := V / V_0, \quad \overline{V_1} := V_1 / (V_1 \cap V_0).
  \end{gather*}
  ($V$ is defined in Subsection \ref{4.1}, while $V_1$ and $V_0$ are defined in Subsection \ref{4.3}.)

  \begin{theorem}
    Recall that $m_{s_1s_2} \ge 4$ and there is a circuit in the Coxeter graph.
    If $m_{s_1s_2} > 4$, then $\overline{V}$ is an irreducible representation of $W$.
    If $m_{s_1s_2} = 4$, then $\overline{V_1}$ is an irreducible representation of $W$.
    Moreover, $\overline{V}$ and $\overline{V_1}$ are both infinite dimensional.
  \end{theorem}

  \begin{proof}
    From Lemma \ref{lem-3.4} we already know that the representations mentioned in this theorem are irreducible.
    Thus, it suffices to show they are both infinite dimensional.
    Obviously, $\dim \overline{V} \ge \dim \overline{V_1}$.
    Thus, we only need to show $\dim \overline{V_1} = \infty$.
    Let $s \in S$ be arbitrary.
    Then $p^{-1} (s) \cap S_1^\prime$ is an infinite set.
    Let
    \begin{equation*}
      U := \bigoplus_{a \in p^{-1} (s) \cap S_1^\prime} \mathbb{C} \alpha_a.
    \end{equation*}
    Then $U \subseteq V_1$, and $\dim U = \infty$.
    For any $0 \neq v \in U$, we have $s \cdot v = -v$.
    Thus, $U \cap V_0 = 0$, and $\overline{V_1} = V_1 / (V_1 \cap V_0)$ is infinite dimensional.
  \end{proof}

  \begin{remark}
    In the theorem, it is necessary to take the quotient by $V_0$ ($V_0$ might be nonzero). For example, consider the universal covering in Figure \ref{fig7},
    where $p(s_i^\prime) = p(s_i^{\prime\prime}) = s_i$, $p(t_i^\prime) = t_i$, $p(r_i^\prime) = r_i$, $\forall i$, and all edges in the Coxeter graph are labelled by $3$ (hence omitted) except the triangle $s_0s_1s_2$.
    In this situation, the vector
    $$\alpha_{t_1^\prime} + \alpha_{t_2^\prime} + \alpha_{t_3^\prime} + \alpha_{t_4^\prime} + 2 \alpha_{t_0^\prime} - 2 \alpha_{r_0^\prime} - \alpha_{r_1^\prime} - \alpha_{r_2^\prime} - \alpha_{r_3^\prime} - \alpha_{r_4^\prime} $$
    is fixed by every $s_i$, $r_i$ and $t_i$.
    Besides, under the setting of Section \ref{sec-3} it is also necessary to take the quotient by $V_0$ in Theorem \ref{thm3.2} (an example can be given in a similar way).
  \end{remark}

  \begin{figure} [ht]
    \centering
    \begin{tikzpicture}
      \node [circle, draw, inner sep=2pt, label=below:$s_1^\prime$] (s1) at (-2,-0.5) {};
      \node [circle, draw, inner sep=2pt, label=below:$s_2^\prime$] (s2) at (2,-0.5) {};
      \node [circle, draw, inner sep=2pt, label=below:$s_0^\prime$] (s0) at (0,1) {};
      \node [circle, draw, inner sep=2pt, label=below:$s_1^{\prime\prime}$] (s1p) at (3.5,-0.5) {};
      \node [circle, draw, inner sep=2pt, label=below:$s_2^{\prime\prime}$] (s2p) at (-3.5,-0.5) {};
      \node [circle, draw, inner sep=2pt, label=left:$t_1^\prime$] (t1) at (-3.2,1) {};
      \node [circle, draw, inner sep=2pt, label=below:$t_0^\prime$] (t0) at (-1.6,1) {};
      \node [circle, draw, inner sep=2pt, label=left:$t_2^\prime$] (t2) at ($(-1.6,1) + (135:1.6)$) {};
      \node [circle, draw, inner sep=2pt, label=left:$t_3^\prime$] (t3) at ($(-1.6,1) + (90:1.6)$) {};
      \node [circle, draw, inner sep=2pt, label=left:$t_4^\prime$] (t4) at ($(-1.6,1) + (45:1.6)$) {};
      \node [circle, draw, inner sep=2pt, label=below:$r_0^\prime$] (r0) at (1.6,1) {};
      \node [circle, draw, inner sep=2pt, label=right:$r_1^\prime$] (r1) at (3.2,1) {};
      \node [circle, draw, inner sep=2pt, label=right:$r_2^\prime$] (r2) at ($(1.6,1) + (45:1.6)$) {};
      \node [circle, draw, inner sep=2pt, label=right:$r_3^\prime$] (r3) at ($(1.6,1) + (90:1.6)$) {};
      \node [circle, draw, inner sep=2pt, label=right:$r_4^\prime$] (r4) at ($(1.6,1) + (135:1.6)$) {};
      \node [label=left:$\cdots$] (left) at (-4.5,-0.5) {};
      \node [label=right:$\cdots$] (right) at (4.5,-0.5) {};
      \graph {(t1) -- (t0) -- (s0) -- (s1) --(s2p) -- (left); (right) -- (s1p) -- (s2) -- (s0) -- (r0) -- (r1); (t2) -- (t0) -- (t3); (t0) -- (t4); (r2) -- (r0) -- (r3); (r0) -- (r4);};

      \draw [->] (0, -1) -- (0, -2.3);
      \node [label=right:$p$] (mapp) at (-0.1, -1.65) {};

      \node [circle, draw, inner sep=2pt, label=left:$s_1$] (s1) at (-2,-6) {};
      \node [circle, draw, inner sep=2pt, label=right:$s_2$] (s2) at (2,-6) {};
      \node [circle, draw, inner sep=2pt, label=below:$s_0$] (s0) at (0,-4.5) {};
      \node [circle, draw, inner sep=2pt, label=left:$t_1$] (t1) at (-3.2,-4.5) {};
      \node [circle, draw, inner sep=2pt, label=below:$t_0$] (t0) at (-1.6,-4.5) {};
      \node [circle, draw, inner sep=2pt, label=left:$t_2$] (t2) at ($(-1.6,-4.5) + (135:1.6)$) {};
      \node [circle, draw, inner sep=2pt, label=left:$t_3$] (t3) at ($(-1.6,-4.5) + (90:1.6)$) {};
      \node [circle, draw, inner sep=2pt, label=left:$t_4$] (t4) at ($(-1.6,-4.5) + (45:1.6)$) {};
      \node [circle, draw, inner sep=2pt, label=below:$r_0$] (r0) at (1.6,-4.5) {};
      \node [circle, draw, inner sep=2pt, label=right:$r_1$] (r1) at (3.2,-4.5) {};
      \node [circle, draw, inner sep=2pt, label=right:$r_2$] (r2) at ($(1.6,-4.5) + (45:1.6)$) {};
      \node [circle, draw, inner sep=2pt, label=right:$r_3$] (r3) at ($(1.6,-4.5) + (90:1.6)$) {};
      \node [circle, draw, inner sep=2pt, label=right:$r_4$] (r4) at ($(1.6,-4.5) + (135:1.6)$) {};
      \graph {(t1) -- (t0) -- (s0) -- (s1) -- (s2) -- (s0) -- (r0) -- (r1); (t2) -- (t0) -- (t3); (t0) -- (t4); (r2) -- (r0) -- (r3); (r0) -- (r4);};
      \node  (m12) at (0,-5.85) {$m_{12}$};
      \node  (m01) at (-0.7,-5.4) {$m_{01}$};
      \node  (m02) at (0.7,-5.4) {$m_{02}$};
      \end{tikzpicture}
    \caption{A universal covering}\label{fig7}
  \end{figure}

\section{Another example} \label{sec-5}
  In this section, we construct an infinite-dimensional irreducible representation of a specific Coxeter group whose Coxeter graph is a tree.
  Suppose the Coxeter graph $G$ of $(W,S)$ is
  
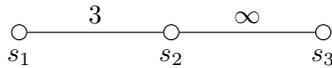
\begin{equation*}
     \begin{tikzpicture}
      \node [circle, draw, inner sep=2pt, label=below:$s_1$] (s1) at (-2,0) {};
      \node [circle, draw, inner sep=2pt, label=below:$s_2$] (s2) at (0,0) {};
      \node [circle, draw, inner sep=2pt, label=below:$s_3$] (s3) at (2,0) {};
      \node [above] (m12) at (-1,0) {$3$};
      \node [above] (m23) at (1,0) {$\infty$};
      \graph{(s1) -- (s2) -- (s3);};
     \end{tikzpicture}
  \end{equation*}
  This Coxeter group is isomorphic to $\operatorname{PGL}(2,\mathbb{Z})$ (see \cite[\S 5.1]{Humphreys90}).

  Recall that the dihedral subgroup $\langle s_1, s_2 \rangle$ has an irreducible representation $\rho_1$ on $\mathbb{C} \beta_{s_1} \oplus \mathbb{C} \beta_{s_2}$  (see Subsection \ref{subsec2.1}).
  There is a basis $\{u,v\}$ of this space such that $u,v$ are eigenvectors of $s_1$ with eigenvalues $+1, -1$, respectively.
  The vectors $u,v$ can be chosen so that $s_2 \cdot u = (3v-u)/2$, $s_2 \cdot v = (u+v)/2$, as illustrated in Figure \ref{fig2}.

  \begin{figure}[ht]
    \centering
    \begin{tikzpicture}
      \draw[dashed] (-1.3,0)--(2,0);
      \draw[dashed] (240:1.3)--(60:2);
      \draw[->] (0,0)--(0,1.5);
      \draw[->] (0,0)--(330:1.5);
      \draw (0.3,0) arc (0:60:0.3);
      \draw[<->] ($(-1.3,0) + (120:0.2)$) arc (120:240:0.2);
      \draw[<->] ($(60:1.9) + (0.2,0)$) arc (0:120:0.2);
      \draw[->] (0,0)--($(0,1.5) + (330:3)$);
      \node[above] (u) at (2.6,0) {$u = \beta_{s_1} + 2 \beta_{s_2}$};
      \node[left] (ar) at (0,1.5) {$v =\beta_{s_1}$};
      \node[below] (at) at (330:1.5) {$\beta_{s_2}$};
      \node[right] (r) at (-2.1,0) {$s_1$};
      \node[right] (t) at (63:2.3) {$s_2$};
      \node[right] (angle) at (0.2,0.3) {$\frac{\uppi}{3}$};
    \end{tikzpicture}
    \caption{The basis vectors $u,v$ of the representation $\rho_1$}\label{fig2}
  \end{figure}

  Let
  \begin{equation*}
      V := \mathbb{C} u_0 \oplus \bigoplus_{i \in \mathbb{N}_{>0}} (\mathbb{C} u_i \oplus \mathbb{C} v_i)
  \end{equation*}
  be the vector space with basis $\{u_i, v_j \mid i \in \mathbb{N}, j \in \mathbb{N}_{>0}\}$, and let $s_1, s_2$ act on $V$ by
  \begin{gather*}
      s_1 \cdot u_0 = s_2 \cdot u_0 = u_0, \\
      s_1 \cdot u_i = u_i, \quad s_1 \cdot v_i = -v_i, \quad \forall i\in \mathbb{N}_{>0}, \\
      s_2 \cdot u_i = \frac{3v_i-u_i}{2}, \quad s_2 \cdot v_i = \frac{u_i+v_i}{2}, \quad \forall i \in \mathbb{N}_{>0}.
  \end{gather*}
  Then as a representation of $\langle s_1, s_2 \rangle$, $V$ is a direct sum of a trivial representation $\mathds{1}$ and infinite many copies of $\rho_1$.

  To make $V$ be a representation of $W$, we only need to find an involution on $V$ commuting with $s_1$, and let $s_3$ act by this involution.
  Let
  \begin{align*}
      s_3 \cdot u_{2k} & = u_{2k+1}, \quad s_3 \cdot u_{2k+1} = u_{2k}, \quad \forall k \in \mathbb{N}, \\
      s_3 \cdot v_{2k} & = v_{2k-1}, \quad s_3 \cdot v_{2k-1} = v_{2k}, \quad \forall k \in \mathbb{N}_{>0}.
  \end{align*}
  Intuitively, $s_3$ permutes these basis vectors:
  \begin{equation*}
      \xymatrix{u_0 \ar@/^/@{<->}[r] & u_1 \ar@{--}[d] & u_2 \ar@{--}[d] \ar@/^/@{<->}[r] & u_3 \ar@{--}[d] & u_4 \ar@{--}[d] \ar@/^/@{<->}[r] & \cdots \\
      & v_1 \ar@/^/@{<->}[r] & v_2 & v_3 \ar@/^/@{<->}[r] & v_4 & \cdots}
  \end{equation*}
  Obviously, the action of $s_3$ is an involution, and commutes with the action of $s_1$.
  Thus, $V$ forms a representation of $W$.

  Similar to the proof of Lemma \ref{lem-3.4}, $u_0$ lies in any nonzero subrepresentation of $V$.
  Note that $u_0$ generates the whole $V$.
  So $V$ is an irreducible representation of $W$.
  \begin{remark}
    Similar constructions also give irreducible representations of infinite dimension when $m_{s_1s_2} = 3$ is replaced by larger integers.
  \end{remark}

\section{Appendix: A sketched proof of Theorem \ref{thm1.1}} \label{sec-app}
This proof is given by an anonymous referee of a previous version of this paper.
The proof uses the following fact.

\begin{lemma}[{\cite[Corollary 2]{MV00}}] \label{lem-app} 
  Suppose $(W,S)$ is an infinite non-affine irreducible Coxeter group of finite rank.
  Then there exists a subgroup $Y^\prime \subseteq W$ of finite index and a surjective homomorphism $\varphi^\prime: Y^\prime \twoheadrightarrow F^\prime$, where $F^\prime$ is a non-abelian free group.
\end{lemma}

Now we can prove Theorem \ref{thm1.1}, and we only need to prove the ``only if'' part.

Suppose $(W,S)$ is infinite and non-affine.
Let $Y^\prime$, $\varphi^\prime$, $F^\prime$ be as in Lemma \ref{lem-app}.
Let $Y$ be the intersection of all conjugates of $Y^\prime$ in $W$.
Then $Y$ is a normal subgroup of $W$ of finite index.
Since $Y$ is of finite index in $Y^\prime$, the image $\varphi^\prime(Y)$ is a subgroup of $F^\prime$ of finite index.
Therefore, $\varphi^\prime(Y)$ is also a non-abelian free group (see, for example,  \cite[Theorem 85.1]{Munkres00}).
Let $X_1, X_2$ be two free generators of $\varphi^\prime(Y)$, and let $F = \langle  X_1, X_2 \rangle$.
Then $F$ is a free group of rank two, and we have a surjection $\varphi^\prime(Y) \twoheadrightarrow F$.
By composing this surjection with the map $\varphi^\prime|_Y$, we obtain a surjective homomorphism $\varphi: Y \twoheadrightarrow F$.

Let $V$ be an infinite-dimensional irreducible representation of $F$ (for example, $V := \bigoplus_{n \in \mathbb{Z}} \mathbb{C} \alpha_{n}$, and let $X_1 \cdot \alpha_n = \alpha_{n+1}$, $X_2 \cdot \alpha_n = 2^n \alpha_{n+1}$, as in the proof of Theorem \ref{thm3.2}).
By pulling back this $F$-representation along $\varphi$, $V$ becomes an irreducible representation of $Y$.

We consider the induced representation $V_Y^W: = \mathbb{C}[W] \otimes_{\mathbb{C}[Y]} V$ of $W$, and we choose $w_1, \dots, w_n$ as coset representatives for $Y$ in $W$.
Then $V_Y^W = \bigoplus_{1 \le i \le n} w_i \otimes V$ as a vector space.
Since $Y$ is a normal subgroup of $W$, each summand $w_i \otimes V$ is an irreducible representation of $Y$.
Thus, $V_Y^W$ is a semisimple Noetherian and Artinian $\mathbb{C}[Y]$-module, and hence a Noetherian and Artinian $\mathbb{C}[W]$-module.
Consequently, there exists an irreducible $W$-representation $M \subseteq V_Y^W$ as a subrepresentation.

If we view $M$ as a $\mathbb{C}[Y]$-module, then, by semi-simplicity of the $\mathbb{C}[Y]$-module $V_Y^W$ and Schur's lemma, $M$ is isomorphic to a direct sum of some $\mathbb{C}[Y]$-modules $w_i \otimes V$.
In particular, $M$ is infinite dimensional.
Theorem \ref{thm1.1} is proved.

\section*{Acknowledgement}
The author would like to thank professor Nanhua Xi for useful discussions. The author is also grateful to the anonymous referee who provided a proof of Theorem \ref{thm1.1}, and to the other referee for useful comments which improved this paper a lot.

\bibliographystyle{amsplain}
\bibliography{inf-dim-rep}

\end{document}